\documentclass[11pt]{article}
\usepackage[T1]{fontenc}
\usepackage[utf8]{inputenc}
\usepackage{lmodern}
\usepackage[margin=1in]{geometry}
\usepackage{amsmath,amssymb,amsfonts,amsthm}
\usepackage[authoryear,round]{natbib}
\usepackage{microtype}
\usepackage{graphicx}
\usepackage{tikz}
\usetikzlibrary{patterns}
\usepackage[hidelinks]{hyperref}

\newtheorem{theorem}{Theorem}
\newtheorem{proposition}{Proposition}
\newtheorem{lemma}{Lemma}

\newcommand{\X}{\mathcal X}
\newcommand{\one}{\mathbf 1}
\newcommand{\dd}{\,\mathrm d}

\title{A Partial-Matching Proof of the Neyman--Pearson Lemma}
\author{Rami V. Tabri\\[3pt]
\normalsize Department of Econometrics and Business Statistics\\
\normalsize Monash University, Melbourne, Australia\\
\normalsize \texttt{rami.tabri@monash.edu}}
\date{}

\begin{document}
\maketitle
\begin{abstract}
Testing seeks to distinguish distributions, whereas matching seeks to identify their common mass. I show that this apparent opposition provides a direct proof of the Neyman--Pearson lemma. Maximal partial matching between the alternative distribution and a scaled null distribution separates them into a common component and two mutually singular residuals. An elementary dual problem then turns this separation into a testing rule: complementary slackness forces rejection on the alternative residual and non-rejection on the scaled-null residual, leaving freedom to randomize at the likelihood-ratio threshold. Thus the optimal test emerges from maximal matching rather than being proposed in advance. The matching profile also determines the optimal power envelope.
\end{abstract}

\noindent\textbf{Keywords:} Neyman--Pearson lemma; partial transport; maximal coupling; complementary slackness; randomized test.
\medskip

\section{Introduction}
\label{sec:intro}
The Neyman--Pearson lemma, established by
\citet{neymanpearson1933}, is one of the foundational
results of mathematical statistics. It characterizes
the most powerful test of a simple null hypothesis
against a simple alternative at a prescribed
significance level, thereby providing a fundamental
justification for likelihood-ratio testing.
The classical proof establishes optimality by comparing
a likelihood-ratio test with an arbitrary competitor.
The result can also be obtained through variational
or Lagrangian arguments applied directly to the
optimization problem over tests.

In this note, I approach the problem from a different
perspective, motivated by the relationship between
hypothesis testing and the common mass of two measures.
A test seeks to distinguish the null and alternative
distributions, while a matching problem seeks to
identify how much of their mass can be placed in
correspondence. The mass that cannot be matched
provides a natural measure of their distinguishability.
This suggests studying hypothesis testing by first
separating the common and unmatched components of
the two distributions, rather than beginning with
a candidate rejection rule.

To connect maximal matching with the trade-off between
testing errors, I compare the alternative distribution
$P_1$ with a scaled null distribution $\lambda P_0$,
where $\lambda\geq 0$. The scaling parameter, $\lambda$, will play
the role of the weight on type-I error in the associated
dual problem. I seek the greatest common submeasure
of these two measures, or equivalently the maximum mass
that can be matched along the diagonal. Removing this
common component leaves two mutually singular residual
measures. The residual associated with the alternative
identifies where rejection is optimal, while the
residual associated with the scaled null identifies
where non-rejection is optimal. Where both residuals
vanish, the decision remains unrestricted, providing
precisely the freedom needed for threshold
randomization.

The connection with hypothesis testing is made
rigorous through an elementary two-potential dual
formulation of the maximal partial-matching problem.
After normalization, one dual potential takes values
in $[0,1]$ and can be interpreted as a randomized
test. Complementary slackness determines its values
on the unmatched residuals, recovering the
likelihood-ratio rule. Choosing the scaling parameter
and the randomization probability to satisfy the
significance constraint then yields the
Neyman--Pearson lemma.

The construction is closely related to the classical
Jordan decomposition of a signed measure: the unmatched
residuals are the positive and negative parts of
$P_1-\lambda P_0$. The identities connecting measure
overlap, testing risk, and convex duality are well
known; see, for example, \citet{cvitanic2001},
\citet{gushchin2020}, and
\citet{mullhaupt2025}. My contribution is
therefore not a new optimality theorem, but a
self-contained proof that makes explicit how maximal
matching, residual separation, and complementary
slackness together determine the Neyman--Pearson test.
\section{A maximal partial-matching problem}
\label{sec:matching}
I begin with a partial transport formulation
that seeks to match as much mass as possible
between the alternative distribution and a
scaled null distribution, without requiring
all available mass to be transported. Let $(\mathcal X,\mathcal A)$ be a standard Borel
space, and let $P_0,P_1$ be probability measures
on it. Write
\[
\Delta=\{(x,x):x\in\mathcal X\}
\]
for the diagonal of $\mathcal X^2$. For
$\lambda\geq0$, define the set of partial
transport plans
\[
\Pi_{\leq}(P_1,\lambda P_0)
=
\bigl\{
\pi\in\mathcal M_+(\mathcal X^2):
\pi_1\leq P_1,\;
\pi_2\leq\lambda P_0
\bigr\},
\]
where $\mathcal M_+(\mathcal X^2)$ denotes
the finite nonnegative measures on
$\mathcal X^2$, and $\pi_1,\pi_2$ are the
marginals of $\pi$. Consider
\[
M(\lambda)
=
\sup_{\pi\in\Pi_{\leq}(P_1,\lambda P_0)}
\int_{\mathcal X^2}
\mathbf 1_\Delta(x,y)\,\pi(dx,dy).
\]
The objective maximizes mass matched at
identical outcomes. Equivalently, it rewards
zero-cost matches under the diagonal cost
$c_\Delta(x,y)=\mathbf 1\{x\ne y\}$.

For any feasible $\pi$, its restriction to
$\Delta$ is also feasible and has the same
objective value. A transport plan concentrated
on $\Delta$ has identical marginals, given
by a measure $\gamma$ satisfying
$\gamma\leq P_1$ and
$\gamma\leq\lambda P_0$. Conversely, every
such $\gamma$ induces a diagonal transport
plan. Thus the partial transport problem
reduces to the maximal common-submeasure
problem
\begin{equation}
\label{eq:primal}
M(\lambda)
=
\sup_{\gamma\in\mathcal M_+(\mathcal X)}
\bigl\{
\gamma(\mathcal X):
0\leq\gamma\leq P_1,\;
\gamma\leq\lambda P_0
\bigr\}.
\end{equation}
Here inequalities between measures mean
domination on every measurable set.

Set $\rho=P_0+P_1$ and
$p_i=dP_i/d\rho$, $i=0,1$.
The greatest common submeasure is explicit:
\begin{equation}\label{eq:overlap}
 \dd\gamma_\lambda^*
 =\min\{p_1,\lambda p_0\}\dd\rho,
 \quad
 M(\lambda)=
 \int\min\{p_1,\lambda p_0\}\dd\rho.
\end{equation}
Indeed, every feasible $\gamma$ is dominated by $\rho$
and has a density bounded by both $p_1$ and
$\lambda p_0$. Consequently,
$\gamma\leq\gamma_\lambda^*$, and the displayed
measure is the unique maximizer. In particular, the primal
optimum is obtained by a greatest element in the measure order.

The two-potential dual can be derived directly
from the partial transport formulation.
Associate nonnegative bounded measurable
potentials $f$ and $g$ with the marginal
domination constraints. The Lagrangian is
\begin{align*}
\mathcal L(\pi;f,g)
=\int f\,dP_1+\lambda\int g\,dP_0
+\int_{\mathcal X^2}
\bigl[\mathbf 1_\Delta(x,y)-f(x)-g(y)\bigr]
\,d\pi(x,y).
\end{align*}
Its supremum over finite nonnegative measures
$\pi$ is finite if and only if
$f(x)+g(y)\geq\mathbf 1_\Delta(x,y)$
for every $(x,y)\in\mathcal X^2$.
Since $f,g\geq0$, this condition reduces to
$f+g\geq1$ on $\mathcal X$. This leads to
the two-potential dual problem
\begin{equation}\label{eq:dual}
 D(\lambda)=
 \inf_{\substack{f,g:\X\to[0,\infty)\\
                 f+g\geq1}}
 \left\{
 \int f\dd P_1+
 \lambda\int g\dd P_0
 \right\},
\end{equation}
where the potentials are measurable. For any feasible
$\gamma$ and $(f,g)$ with finite dual objective,
\begin{align}\label{eq:gap}
 \int f\dd P_1+
 \lambda\int g\dd P_0-\gamma(\X)
 =
 \int f\dd(P_1-\gamma)
 +\int g\dd(\lambda P_0-\gamma)+
 \int(f+g-1)\dd\gamma.
\end{align}
Every term on the right-hand side is nonnegative,
establishing weak duality.

\begin{proposition}[Matching duality]
\label{prop:duality}
For every $\lambda\geq0$, the primal and dual
problems \eqref{eq:primal} and \eqref{eq:dual}
have equal values, and both optima are attained:
\(
 D(\lambda)=M(\lambda).
\)
Moreover, the dual problem can be restricted
without loss of optimality to normalized pairs
\[
 (f,g)=(1-g,g),\qquad 0\leq g\leq1.
\]
\end{proposition}

\begin{proof}
For any feasible dual pair $(f,g)$, define
\[
 \widetilde f=\min\{f,1\},
 \qquad
 \widetilde g=1-\widetilde f.
\]
The constraints $f,g\geq0$ and $f+g\geq1$
imply
\[
 0\leq\widetilde f\leq f,\qquad
 0\leq\widetilde g\leq g,\qquad
 \widetilde f+\widetilde g=1.
\]
Thus normalization preserves feasibility and
does not increase the dual objective.

For a normalized pair $(1-g,g)$, the objective
in \eqref{eq:dual} becomes
\begin{equation}\label{eq:risk}
 \int(1-g)\dd P_1+\lambda\int g\dd P_0
 =
 1-\int g(p_1-\lambda p_0)\dd\rho.
\end{equation}
Since $0\leq g\leq1$,
\[
 \int g(p_1-\lambda p_0)\dd\rho
 \leq
 \int(p_1-\lambda p_0)_+\dd\rho.
\]
Equality is attained by
$g=\one_{\{p_1>\lambda p_0\}}$.
Consequently,
\begin{align*}
 D(\lambda)
 =1-\int(p_1-\lambda p_0)_+\dd\rho
 =\int\min\{p_1,\lambda p_0\}\dd\rho
 =M(\lambda).
\end{align*}
Together with \eqref{eq:overlap}, this establishes
strong duality and attainment constructively, without invoking
a general infinite-dimensional duality theorem.
\end{proof}

The matching problem leaves two residual measures:
\begin{equation}\label{eq:residuals}
 R_{1,\lambda}=P_1-\gamma_\lambda^*,
 \qquad
 R_{0,\lambda}=\lambda P_0-\gamma_\lambda^*.
\end{equation}
The next result shows that their separation
characterizes maximal matching and determines
the normalized optimal dual potentials.

\begin{lemma}[Residual separation and optimal decisions]
\label{lem:residual}
Fix $\lambda\geq0$. For any feasible common
submeasure $\gamma$ in \eqref{eq:primal},
define
\[
 R_1^\gamma=P_1-\gamma,\qquad
 R_0^\gamma=\lambda P_0-\gamma.
\]
Then $\gamma=\gamma_\lambda^*$ if and only if
\(
 R_1^\gamma\perp R_0^\gamma.
\)
Furthermore, a measurable function
$g:\X\to[0,1]$ is an optimal normalized
dual potential for \eqref{eq:dual} if and
only if
\begin{equation}\label{eq:cs}
 g=1\quad R_{1,\lambda}\text{-a.e.},
 \qquad
 g=0\quad R_{0,\lambda}\text{-a.e.}
\end{equation}
Equivalently,
\begin{equation}\label{eq:optimal-g}
 g=
 \begin{cases}
 1,&p_1>\lambda p_0,\\
 0,&p_1<\lambda p_0,\\
 \eta,&p_1=\lambda p_0,
 \end{cases}
 \qquad 0\leq\eta\leq1,
\end{equation}
$\rho$-almost everywhere, where $\eta$ may
be any measurable function on the equality set.
\end{lemma}

\begin{proof}
Write $\gamma=h\rho$. Feasibility implies
\[
 0\leq h\leq\min\{p_1,\lambda p_0\}
 \qquad\rho\text{-a.e.}
\]
The residual measures $R_1^\gamma$ and
$R_0^\gamma$ have densities $p_1-h$ and
$\lambda p_0-h$, respectively.

If $\gamma=\gamma_\lambda^*$, then
$h=\min\{p_1,\lambda p_0\}$, and the
residual densities become
\[
 (p_1-\lambda p_0)_+,\qquad
 (\lambda p_0-p_1)_+.
\]
Their positivity sets are disjoint, so
$R_1^\gamma\perp R_0^\gamma$.

Conversely, suppose
$R_1^\gamma\perp R_0^\gamma$.
Because both residual measures are
absolutely continuous with respect to
$\rho$, their densities cannot both
be strictly positive on a set of
positive $\rho$-measure. Hence
\[
 (p_1-h)(\lambda p_0-h)=0
 \qquad\rho\text{-a.e.}
\]
Since $h$ does not exceed either
$p_1$ or $\lambda p_0$, this implies
\[
 h=\min\{p_1,\lambda p_0\}
 \qquad\rho\text{-a.e.}
\]
Thus $\gamma=\gamma_\lambda^*$.

Now let $(1-g,g)$ be a normalized
feasible dual pair. Substituting
$\gamma_\lambda^*$ into
\eqref{eq:gap} gives
\begin{align*}
 \int(1-g)\dd P_1+
 \lambda\int g\dd P_0-M(\lambda)
 =
 \int(1-g)\dd R_{1,\lambda}
 +\int g\dd R_{0,\lambda}.
\end{align*}
By Proposition~\ref{prop:duality},
strong duality holds. Both terms on
the right are nonnegative, so $g$
is optimal if and only if both
vanish. This is equivalent to
\eqref{eq:cs}.

Finally, $R_{1,\lambda}$ has
strictly positive density on
$\{p_1>\lambda p_0\}$, whereas
$R_{0,\lambda}$ has strictly
positive density on
$\{p_1<\lambda p_0\}$. Thus
\eqref{eq:cs} forces $g=1$
and $g=0$, respectively, on
these regions, up to $\rho$-null
sets. On $\{p_1=\lambda p_0\}$,
both residual measures vanish,
so any measurable value in
$[0,1]$ is permitted.
This proves \eqref{eq:optimal-g}.
\end{proof}

The residual measures in \eqref{eq:residuals} are the positive
and negative parts in the Jordan decomposition
\[
 P_1-\lambda P_0=R_{1,\lambda}-R_{0,\lambda}.
\]
The lemma therefore reads the familiar sign decomposition as a
property of maximal matching: once all common mass has been removed,
the unmatched parts are mutually singular. Complementary slackness
then fixes the dual decision on those parts. The equality set
$\{p_1=\lambda p_0\}$ is where both residual densities vanish;
it need not coincide with the support of the matched measure.

\section{The Neyman--Pearson lemma from matching duality}
\label{sec:np}

A randomized test is a measurable function
$\varphi:\X\to[0,1]$. Its size and power are
\[
 P_0\varphi=\int\varphi\dd P_0,
 \qquad
 P_1\varphi=\int\varphi\dd P_1.
\]
For $\alpha\in[0,1]$, define the optimal
power function
\begin{equation}\label{eq:power}
 \beta(\alpha)=
 \sup_{\substack{0\leq\varphi\leq1\\
                 P_0\varphi\leq\alpha}}
 P_1\varphi.
\end{equation}

The maximal partial-matching problem supplies a bound on the
power of every test. For $\lambda\geq0$, the decompositions
\[
 P_1=\gamma_\lambda^*+R_{1,\lambda},\qquad
 \lambda P_0=\gamma_\lambda^*+R_{0,\lambda}
\]
imply
\begin{equation}\label{eq:residual-identity}
 P_1\varphi-\lambda P_0\varphi
 =\int\varphi\dd R_{1,\lambda}
 -\int\varphi\dd R_{0,\lambda}.
\end{equation}
Because $0\leq\varphi\leq1$ and both residual
measures are nonnegative,
\[
\int\varphi\,dR_{1,\lambda}
\leq R_{1,\lambda}(\mathcal X)
=1-M(\lambda),
\qquad
\int\varphi\,dR_{0,\lambda}\geq0.
\]
Combining these inequalities with
\eqref{eq:residual-identity} gives
\begin{equation}\label{eq:matching-bound}
P_1\varphi
\leq\lambda P_0\varphi+1-M(\lambda).
\end{equation}
More precisely, for any $\alpha\geq P_0\varphi$,
\begin{align}\label{eq:power-gap}
 \lambda\alpha+1-M(\lambda)-P_1\varphi
 =\lambda(\alpha-P_0\varphi)
 +\int(1-\varphi)\dd R_{1,\lambda}
+\int\varphi\dd R_{0,\lambda}.
\end{align}

All three terms on the right-hand side of \eqref{eq:power-gap} are
nonnegative. Thus the decomposition establishes not
only a bound on the power of every feasible test, but
also the conditions under which that bound is attained.
By Lemma~1, the last two terms vanish exactly when
$\varphi$ is an optimal normalized matching-dual
potential: rejection is required on the unmatched
alternative residual, non-rejection on the unmatched
scaled-null residual, and randomization is permitted
where both residuals vanish. For $\lambda>0$, equality
in (12) additionally requires $P_0\varphi=\alpha$;
when $\lambda=0$, the size constraint need not bind.
Consequently, optimality for the weighted testing
problem at a fixed $\lambda$ and optimality under a
prescribed significance constraint are connected
through the equality conditions in \eqref{eq:power-gap}.

\begin{theorem}[Neyman--Pearson lemma]
\label{thm:np}
For every $\alpha\in[0,1]$, a most
powerful test of size at most $\alpha$
exists. For $0<\alpha<1$, there are
$\lambda_\alpha\geq0$ and
$\eta_\alpha\in[0,1]$ such that
\begin{equation}\label{eq:np}
 \varphi_\alpha=
 \one_{\{p_1>\lambda_\alpha p_0\}}
 +\eta_\alpha
 \one_{\{p_1=\lambda_\alpha p_0\}}
\end{equation}
has size $P_0\varphi_\alpha=\alpha$
and is most powerful among all tests
of size at most $\alpha$.

Moreover, $\varphi_\alpha$ is a
normalized optimal dual potential
for \eqref{eq:dual} at
$\lambda_\alpha$, and
\begin{equation}\label{eq:envelope}
 \beta(\alpha)=
 \inf_{\lambda\geq0}
 \left\{
 \lambda\alpha+1-M(\lambda)
 \right\},
 \qquad 0\leq\alpha\leq1.
\end{equation}
\end{theorem}

\begin{proof}
Fix $0<\alpha<1$. On $\{p_0>0\}$,
define the likelihood ratio $L=p_1/p_0$. Since $L$ is
finite $P_0$-almost surely, the
threshold
\begin{equation}\label{eq:threshold}
 \lambda_\alpha=
 \inf\left\{
 t\geq0:P_0(L>t)\leq\alpha
 \right\}
\end{equation}
is finite. The likelihood ratio is used only to select the
threshold under $P_0$. The rejection region itself
is defined by the density comparison
$p_1>\lambda_\alpha p_0$, which also includes
$\{p_0=0,\ p_1>0\}$. The monotonicity and
continuity properties of the survival
function imply
\begin{equation}\label{eq:quantile}
 P_0(L>\lambda_\alpha)
 \leq\alpha
 \leq P_0(L\geq\lambda_\alpha).
\end{equation}

Let
\(
 A_+=\{p_1>\lambda_\alpha p_0\},
 \;
 A_0=\{p_1=\lambda_\alpha p_0\}.
\)
Since $P_0(\{p_0=0\})=0$,
\eqref{eq:quantile} gives
\[
 P_0(A_+)\leq\alpha
 \leq P_0(A_+)+P_0(A_0).
\]
Choose
\begin{equation}\label{eq:eta}
 \eta_\alpha=
 \begin{cases}
 \displaystyle
 \frac{\alpha-P_0(A_+)}{P_0(A_0)},
 &P_0(A_0)>0,\\[8pt]
 0,&P_0(A_0)=0.
 \end{cases}
\end{equation}
Then $\eta_\alpha\in[0,1]$, and
the test \eqref{eq:np} satisfies
$P_0\varphi_\alpha=\alpha$.

By Lemma~\ref{lem:residual},
$\varphi_\alpha$ satisfies
\[
 \int(1-\varphi_\alpha)
       \dd R_{1,\lambda_\alpha}=0,
 \qquad
 \int\varphi_\alpha
       \dd R_{0,\lambda_\alpha}=0.
\]
Since $P_0\varphi_\alpha=\alpha$, all three terms in
\eqref{eq:power-gap} vanish. Consequently,
\[
 P_1\varphi_\alpha
 =
 \lambda_\alpha\alpha+
 1-M(\lambda_\alpha).
\]
On the other hand, every test
$\varphi$ with $P_0\varphi\leq\alpha$
satisfies, by \eqref{eq:power-gap},
\[
 P_1\varphi
 \leq
 \lambda_\alpha\alpha+
 1-M(\lambda_\alpha).
\]
Thus $\varphi_\alpha$ is most powerful.

Moreover, since \eqref{eq:power-gap}
holds for every $\lambda\geq0$,
\[
 \beta(\alpha)
 \leq
 \inf_{\lambda\geq0}
 \{\lambda\alpha+1-M(\lambda)\}.
\]
Equality follows from the preceding
construction at $\lambda_\alpha$,
proving \eqref{eq:envelope} for
$0<\alpha<1$.

At $\alpha=0$, every feasible test
vanishes $P_0$-almost surely.
The test
$\varphi_0=\one_{\{p_0=0\}}$
therefore attains
\[
 \beta(0)=P_1(\{p_0=0\}).
\]
Because
\(
 1-M(\lambda)
 =
 \int(p_1-\lambda p_0)_+\dd\rho,
\)
dominated convergence gives
\[
 \lim_{\lambda\to\infty}
 \{1-M(\lambda)\}
 =
 P_1(\{p_0=0\}),
\]
establishing \eqref{eq:envelope}
at $\alpha=0$.

Finally, at $\alpha=1$,
$\varphi_1\equiv1$ has power one.
Because $M(0)=0$, the infimum in
\eqref{eq:envelope} is also one.
\end{proof}

Thus the maximal matching profile $M(\lambda)$,
as the scaling parameter varies, determines the
entire Neyman--Pearson power envelope through
the infimum representation \eqref{eq:envelope}.

The proof begins with the matched and unmatched components of
the two measures. The common component cancels in
\eqref{eq:residual-identity}, while the nonnegative residual terms
in \eqref{eq:power-gap} identify precisely when a test attains
the power bound. The likelihood-ratio rule follows by expressing
these complementary-slackness conditions in terms of densities.

\section{An illustration with threshold randomization}

Take a three-point sample space $\X=\{a,b,c\}$,
with
\[
\begin{array}{c|ccc}
 &a&b&c\\ \hline
 P_0&0.5&0.3&0.2\\
 P_1&0.2&0.3&0.5\\
 P_1/P_0&0.4&1&2.5
\end{array}
\]
Here the dominating measure introduced in Section~2
is $\rho=P_0+P_1$, which assigns masses
$(0.7,0.6,0.7)$ to $(a,b,c)$. Thus
$p_i(x)=P_i(\{x\})/\rho(\{x\})$, $i=0,1$,
and equation~(2) gives
\begin{align*}
\gamma_\lambda^*(\{x\})
=&\rho(\{x\})\min\{p_1(x),\lambda p_0(x)\}\\
=&\min\{P_1(\{x\}),\lambda P_0(\{x\})\}.
\end{align*}
At $\lambda=1$, the optimal matched mass at
$(a,b,c)$ is therefore $(0.2,0.3,0.2)$,
so $M(1)=0.7$.
The unmatched alternative residual has mass $0.3$
at $c$, while the unmatched null residual has
mass $0.3$ at $a$. Both residuals vanish at $b$.

Complementary slackness therefore requires
$g(c)=1$ and $g(a)=0$, while leaving
$g(b)\in[0,1]$ free. For a size
$\alpha=0.35$, the choice $g(b)=1/2$ gives
\[
 P_0g=0.2+0.3(1/2)=0.35,\;
 P_1g=0.5+0.3(1/2)=0.65.
\]
The duality bound is attained:
\[
 \beta(0.35)=0.35+1-M(1)=0.65.
\]
The randomization reflects the freedom of the optimal
matching-dual potential at $b$, where the two scaled
masses agree and both residual measures vanish.

\begin{figure}[htbp]
\centering
\begin{tikzpicture}[
    x=3.4cm,
    y=3.4cm,
    >=stealth,
    font=\footnotesize
]

\fill[gray!15]
    (0,1)
    -- (0.5,0.8)
    -- (0.8,0.5)
    -- (1,0)
    -- (0.5,0.2)
    -- (0.2,0.5)
    -- cycle;

\draw[thin, gray!75]
    (0,1)
    -- (0.5,0.8)
    -- (0.8,0.5)
    -- (1,0)
    -- (0.5,0.2)
    -- (0.2,0.5)
    -- cycle;

\draw[thin, dashed]
    (0,0.7) -- (0.7,0);

\draw[thick]
    (0,1)
    -- (0.2,0.5)
    -- (0.5,0.2)
    -- (1,0);

\draw[line width=1.5pt]
    (0.2,0.5) -- (0.5,0.2);

\draw[->, thick]
    (0,0) -- (1.07,0);

\draw[->, thick]
    (0,0) -- (0,1.07);

\foreach \t in {0.2,0.4,0.6,0.8,1.0} {
    \draw (\t,0.008) -- (\t,-0.008)
        node[below, font=\scriptsize] {\t};

    \draw (0.008,\t) -- (-0.008,\t)
        node[left, font=\scriptsize] {\t};
}

\node[below left, font=\scriptsize]
    at (0,0) {$0$};

\node[font=\footnotesize]
    at (0.53,-0.16) {$x=P_0\varphi$};

\node[rotate=90, font=\footnotesize]
    at (-0.20,0.53) {$z=1-P_1\varphi$};

\node[font=\scriptsize, anchor=south]
    at (0.53,1.10)
    {$x+z=M(1)=0.7$};

\fill (0.2,0.5) circle (0.8pt);
\fill (0.5,0.2) circle (0.8pt);

\fill (0.35,0.35) circle (1.5pt);

\draw[dotted, thin]
    (0.35,0) -- (0.35,0.35)
    -- (0,0.35);

\node[above right, font=\scriptsize,
      inner sep=2pt]
    at (0.35,0.35)
    {$(0.35,0.35)$};

\end{tikzpicture}

\caption{Attainable risk region for the three-point
example. The solid lower boundary is the efficient
risk frontier. At $\lambda=1$, the dashed supporting
line $x+z=M(1)=0.7$ exposes the segment joining
$(0.2,0.5)$ and $(0.5,0.2)$, corresponding to
threshold randomization at $b$. The marked point
$(0.35,0.35)$ represents the most powerful test
of size $0.35$.}
\label{fig:matching-risk}
\end{figure}
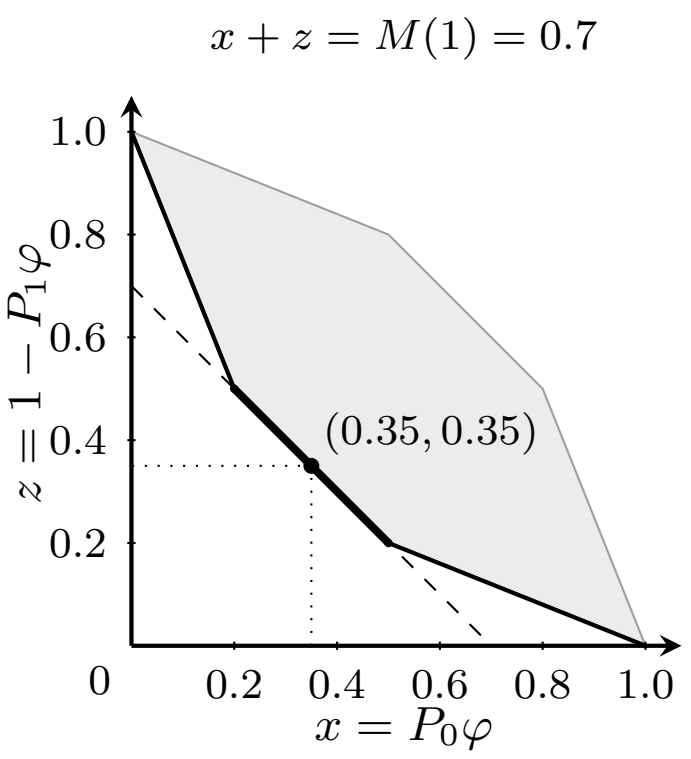

Figure~\ref{fig:matching-risk} illustrates
the geometric interpretation of maximal
partial matching. 

I now illustrate the geometric connection
between maximal matching and the trade-off
between testing errors. For a randomized
test $\varphi$, the probabilities
$P_0\varphi$ and $1-P_1\varphi$ are its
type-I and type-II errors, respectively.
By Proposition~\ref{prop:duality}, the
normalized matching-dual problem can be
written as
\[
M(\lambda)
=
\inf_{0\leq\varphi\leq1}
\bigl\{
\lambda P_0\varphi+1-P_1\varphi
\bigr\}.
\]
Thus maximal matching determines the
minimum weighted sum of the two testing
errors, with $\lambda$ controlling their
relative importance. In particular,
at $\lambda=1$ the objective is the
sum of the type-I and type-II error
probabilities.

To visualize this connection, let
\(
\mathcal R
=
\left\{
(P_0\varphi,1-P_1\varphi):
0\leq\varphi\leq1
\right\}
\)
denote the attainable risk region. Each
point $(x,z)\in\mathcal R$ represents the
two error probabilities associated with
some randomized test. The region is
convex, since randomization between
tests produces convex combinations
of their error probabilities.

The preceding dual representation gives
\[
M(\lambda)
=
\inf_{(x,z)\in\mathcal R}
\{\lambda x+z\}.
\]
Geometrically, this means that
$\mathcal R$ lies on or above the line
\[
z=M(\lambda)-\lambda x,
\]
with equality attained at one or more
points of the region. Hence this line
supports $\mathcal R$ from below.
Its vertical intercept is the maximal
matched mass $M(\lambda)$, and its
slope is $-\lambda$.

As $\lambda$ varies, the supporting
lines identify different trade-offs
between the two error probabilities.
Their points of contact correspond
precisely to tests minimizing the
weighted testing risk, or equivalently,
to optimal normalized dual potentials
in Proposition~\ref{prop:duality}.
By Lemma~\ref{lem:residual}, these
tests are characterized by complementary
slackness on the unmatched residual
measures.

In the present example, $\lambda=1$
gives $M(1)=0.7$. The supporting line
$x+z=0.7$ touches the efficient
boundary along a segment, reflecting
the freedom to randomize at $b$,
where the two scaled distributions
have equal mass. The choice
$\varphi(b)=1/2$ yields the point
$(0.35,0.35)$ on this segment.

\section{Discussion}

The matching formulation reverses the usual order of
reasoning in the Neyman--Pearson argument. Rather
than beginning with a likelihood-ratio rejection
region and verifying its optimality, I first
identify the greatest common submeasure of the
alternative and scaled null distributions.
The resulting residual separation determines the
optimal dual decisions through complementary
slackness. Thus the likelihood-ratio rule emerges
from maximal matching rather than serving as the
starting point of the proof.

The construction is closely related to maximal
coupling. At $\lambda=1$, the common mass
$M(1)$ equals the largest possible probability
that two random variables with distributions
$P_0$ and $P_1$ agree under a coupling.
Writing
$\|P_1-P_0\|_{\mathrm{TV}}
=\sup_{A\in\mathcal A}|P_1(A)-P_0(A)|$
for total variation distance, it follows that
\[
M(1)=1-\|P_1-P_0\|_{\mathrm{TV}}.
\]
For general $\lambda\geq 0$, the matching
problem extends this overlap interpretation
to the finite measures $P_1$ and $\lambda P_0$,
without requiring equal total masses.

The overlap profile also admits an interpretation
in terms of the positive-part, or hockey-stick,
divergence. Recall that maximal matching leaves
an alternative residual $R_{1,\lambda}$ with
density $(p_1-\lambda p_0)_+$. Its total mass
is therefore $1-M(\lambda)$. On the other hand,
for any measurable $g:\mathcal X\to[0,1]$,
\[
P_1g-\lambda P_0g
=
\int g(p_1-\lambda p_0)\,d\rho.
\]
This expression is maximized by taking $g=1$
where $p_1>\lambda p_0$ and $g=0$ where
$p_1<\lambda p_0$. Consequently,
\begin{equation}
\label{eq:hockey_stick_discussion}
1-M(\lambda)
=
\int (p_1-\lambda p_0)_+\,d\rho
=
\sup_{0\leq g\leq1}
\{P_1g-\lambda P_0g\}.
\end{equation}
Thus the unmatched alternative mass is also
the greatest excess of power over $\lambda$
times size attainable by a randomized test.

This interpretation helps explain the
power-envelope representation
\eqref{eq:envelope}. For each $\lambda\geq0$,
the preceding identity gives an upper bound
$\lambda\alpha+1-M(\lambda)$ on the power
of any test of size at most $\alpha$.
The Neyman--Pearson construction shows that
the infimum of these bounds is attained
for $0<\alpha<1$, thereby recovering the
optimal power function. The underlying
variational identities are classical; the
distinctive feature of the present approach
is that they arise from maximal matching
and the separation of its residual measures.

The matching interpretation also provides
a direct explanation of the data-processing
inequality. Let $K$ be a Markov kernel
from $\mathcal X$ to another standard Borel
space $\mathcal X'$, representing a possibly
randomized transformation of the observation.
If $\gamma\leq P_1,\lambda P_0$ is a feasible
common submeasure, then
$\gamma K\leq P_1K,\lambda P_0K$ and
$(\gamma K)(\mathcal X')=\gamma(\mathcal X)$.
Thus any mass that can be matched before
processing can still be matched afterward,
although additional matching may become
possible. Consequently,
\[
M_{P_0K,P_1K}(\lambda)
\geq M_{P_0,P_1}(\lambda).
\]

By \eqref{eq:hockey_stick_discussion},
this is precisely the data-processing
inequality for the hockey-stick divergence:
processing cannot increase
$1-M(\lambda)$. In particular, at
$\lambda=1$ it yields the familiar
contraction of total variation distance,
\[
\|P_1K-P_0K\|_{\mathrm{TV}}
\leq
\|P_1-P_0\|_{\mathrm{TV}}.
\]
More generally, the power-envelope
representation \eqref{eq:envelope} gives
\(
\beta_{P_0K,P_1K}(\alpha)
\leq
\beta_{P_0,P_1}(\alpha).
\)
Hence processing the observations cannot
increase the maximum attainable power at
a given significance level. The matching
perspective unifies these classical facts:
processing cannot decrease common mass,
and therefore cannot improve statistical
distinguishability.

It is useful to clarify the role of transport
in this construction. The diagonal matching
reward introduced in Section~2 records whether
two outcomes agree, rather than measuring
the distance between them. Because unmatched
mass is permitted, transporting mass off the
diagonal provides no benefit, and the problem
reduces to identifying the greatest common
submeasure. In this respect, the construction
differs from Wasserstein transport, where the
cost of moving mass between locations is central.

The resulting proof offers a different way
to understand the Neyman--Pearson lemma.
Maximal matching first separates the
distributions into a common component and
two mutually singular residual measures.
The dual problem then translates this
separation into statistical decisions:
complementary slackness requires rejection
on the unmatched alternative residual and
non-rejection on the unmatched scaled-null
residual, while allowing randomization where
both residuals vanish. Choosing the scaling
parameter to meet the significance constraint
completes the argument. Thus the
likelihood-ratio rule is not imposed at the
outset; it emerges from the structure of
maximal matching and its dual.

\end{document}